\documentclass[a4paper,12pt,reqno]{amsart}

\usepackage{a4}
\usepackage{verbatim}

\usepackage{color}

\usepackage{amssymb}
\usepackage{color}
\usepackage{graphicx}
\usepackage{floatflt}
\usepackage{float} 		

\usepackage{amsthm}
\usepackage{amscd}
\usepackage{hyperref}

\newtheorem{definition}{Definition}
\newtheorem{example}[definition]{Example}

\newtheorem{theorem}[definition]{Theorem}

\newtheorem{remark}[definition]{Remark}

\newtheorem{proposition}[definition]{Proposition}

\def\XXint#1#2#3{{\setbox0=\hbox{$#1{#2#3}{\int}$}
		\vcenter{\hbox{$#2#3$}}\kern-.5\wd0}}

\makeatletter
\@addtoreset{definition}{section}
\@addtoreset{equation}{section}
\makeatother

\DeclareMathOperator{\Res}{Res}


\allowdisplaybreaks[4]

\begin{document}

\title{Genus Permutations and Genus Partitions}

\author{Alexander Hock }

\address{Mathematical Institute, University of Oxford, Andrew Wiles Building, Woodstock Road,
	OX2 6GG, Oxford, UK 
{\itshape email address:} \normalfont  
\texttt{alexander.hock@maths.ox.ac.uk}}

\maketitle

\begin{abstract}
    For a given permutation or set partition there is a natural way to assign a genus. Counting all permutations or partitions of a fixed genus according to cycle lengths or block sizes, respectively, is the main content of this article. After a variable transformation, the generating series are rational functions with poles located at the ramification points in the new variable. The generating series for any genus is given explicitly for permutations and up to genus 2 for set partitions. Extending the topological structure not just by the genus but also by adding more boundaries, we derive the generating series of non-crossing partitions on the cylinder from known results of non-crossing permutations on the cylinder. Most, but not all, outcomes of this article are special cases of already known results, however they are not represented in this way in the literature, which however seems to be the canonical way. To make the article as accessible as possible, we avoid going into details into the explicit connections to Topological Recursion and Free Probability Theory, where the original motivation came from.
\end{abstract}

\markboth{\hfill\textsc\shortauthors}{\textsc{Genus Permutations and Genus Partitions}\hfill}

\section{Introduction}
Permutations and set partitions are very important objects in mathematics and appear in almost all areas. The number of permutations is given by factorials and the number number of partitions by Bell numbers. The ordinary generating series of factorials or Bell numbers only exists as a formal power series, i.e. it has radius of convergence zero. 

One can associate to any permutation or partition a so-called genus, which is an positive integer bounded by the number of elements. Collecting permutations or partitions of a fixed genus according to the lengths of the cycles of the permutations or according to the sizes of the blocks of the partitions is the subject of this article. For each genus, the generating series has a non vanishing radius of convergence. Resumming the non-convergent genus-series is a nontrivial task and highly related to Borel resummation or resurgence, which will be not addressed here. 

For the special case of genus 0, so-called non-crossing permutations and non-crossing partitions coincide, and are counted by Catalan numbers. For higher genus and after a variable transformation, the generating series of genus permutations are rational function with poles just located at the ramification point in the new variable. Theorem \ref{ThmPerm} gives the explicit formula for all genus encoded by integer partitions as a primary example of the work of \cite{Borot:2021thu} but never worked out in such details as in this article. Generalisations in several directions are discussed as well. 

For genus partitions, much less is known. Considering the same variable transformation for genus partitions, the generating series up to genus 2 is given explicitly (based on the work \cite{zuber2023counting}). We find that it has also poles just located at the ramification points in the new variable. Using this rationality, we prove some minor conjectures about genus partitions made in \cite{coquereaux2023counting}. Comparing further permutations and partitions in a more general setting, we provide an explicit formula for non-crossing partitions on the cylinder, which is deduced from the known result of non-crossing permutations on the cylinder.  

The subject of this article is highly connected to the theory of Topological Recursion \cite{Eynard:2007kz} and the theory of Free Probability \cite{Voiculescu1986AdditionOC} and the recently revealed connection between them. We will mention this connection in some remarks throughout the article, but refer the reader to the provided references.

\subsection*{Acknowledgement}
I am grateful to Jean-Bernard Zuber for motivating me to write this article. I want to thank  Robert Coquereaux, Gábor Hetyei, Luca Lionni, James Mingo and Roland Speicher for helpful discussions and literature suggestions. This work was supported through
the Walter-Benjamin fellowship\footnote{``Funded by
	the Deutsche Forschungsgemeinschaft (DFG, German Research
	Foundation) -- Project-ID 465029630}.

\section{Preliminaries}
Taking the set $\{1,...,n\}$, a \textit{permutation} is a bijection from this set to itself, a rearrangement of $\{1,...,n\}$. The set of all permutations on $\{1,...,n\}$ form a group which is called the symmetric group $\mathcal{S}_n$, which has $n!$ elements, i.e. there are $n!$ permutations of $\{1,...,n\}$. A permutation can be decomposed into disjoint cycles, which are the different orbits of $\sigma\in \mathcal{S}_n$. We will call a cycle of length $i$ a $i$-cycle. 

The number of permutations of $\{1,...,n\}$ with $k$ cycles is given by the \textit{Stirling numbers of the first kind} $s(n,k)$ (up the a global sign $(-1)^{n-k}$). The definition is given by the coefficient of the falling factorials
\begin{align*}
    \sum_{k=0}^ns(n,k)x^k=x(x-1)...(x-n+1).
\end{align*}
Let us refine this a bit further. Let $[a]\vdash n$ denote a \textit{partition of the integer} $n$, i.e. $[a]=(a_1,a_2,...,a_k)$ with $\sum_ia_i=n$, $a_i>0$ and $a_i\leq a_{i+1}$. Then, the number of permutations of $\{1,...,n\}$ with cycle type $[a]\vdash n$ is given by
\begin{align*}
    D_{n,[a]}=\frac{n!}{\mathrm{sym}(a)\prod_ia_i},
\end{align*}
where $\mathrm{sym}(a)=\prod_{j=1}^nk_j!$ and $k_j$ is the number of all $a_l=j$ in the integer partition $[a]\vdash n$. The Stirling numbers of the first kind or the factorials are computed from $D_{n,[a]}$ via
\begin{align*}
    (-1)^{n-k}s(n,k)=\sum_{\substack{[a]\vdash n\\ |[a]|=k}}D_{n,[a]},\qquad n!=\sum_{\substack{[a]\vdash n}}D_{n,[a]}.
\end{align*}
Now let $\kappa_i$ be indeterminates, we can define the following moment
\begin{align*}
    \alpha_n=\sum_{[a]\vdash n}D_{n,[a]}\prod_i\kappa_{a_i},
\end{align*}
which associates $\kappa_i$ to a cycle of length $i$ and sums over all integer partitions $[a]\vdash n$. More canonically, we define the generating series $X(y)$ of $\kappa_i$ and the generating series of permutations $W_{per}(x)$
\begin{align}\label{GenKappa}
    X(y)=&\frac{1}{y}+\sum_{i=1}^\infty \kappa_iy^{i-1},\\
    W_{per}(x)=&\frac{1}{x}+\sum_{n=1}\frac{\alpha_n}{x^{n+1}}=\frac{1}{x}+\sum_{n=1}\frac{\sum_{[a]\vdash n}D_{n,[a]}\prod_i\kappa_{a_i}}{x^{n+1}}
\end{align}
and ask what is the relation between the two generating series $X(y)$ and $W_{per}(x)$. 

We should keep three important examples in mind:
\begin{align*}
    &\text{Factorials:} &&\kappa_i=1,\quad  X(y)=\frac{1}{y}+\frac{1}{1-y},\quad W_{per}(x)=\sum_{n=0}^\infty \frac{n!}{x^{n+1}}\\
    &\text{Stirling 1st:} &&\kappa_i=\kappa,\quad X(y)=\frac{1}{y}+\frac{\kappa}{1-y},\quad W_{per}(x)=\sum_{n=0}^\infty \frac{\sum_{k=0}^n(-1)^{n-k}s(n,k) \kappa^k}{x^{n+1}}\\
    &\text{Harer-Zagier:} &&\kappa_i=\delta_{i,2},\quad X(y)=\frac{1}{y}+y,\quad\quad W_{per}(x)=\sum_{n=0}^\infty \frac{(2n-1)!!}{x^{2n+1}}.
\end{align*}
The last example is called Harer-Zagier due to  its relation to the Euler characteristic of the moduli space of curves \cite{MR848681}. The detailed connection will be apparent later.

For any choices of $\kappa_i$, it is obvious that $W_{per}(x)$ can just be understood as a formal power series in $\frac{1}{x}$, which does not converge.\\

Let us now turn to partitions. Taking the set $\{1,...,n\}$, a \textit{set partition} is set of non-empty subsets (blocks) of $\{1,...,n\}$, where each element is in exactly in one of these subsets. Let $P(n)$ be the set of all set partitions of $\{1,...,n\}$. We will write a partition as
$\lambda=\{\lambda_1,...,\lambda_l\}$, where the $\lambda_i$'s are the blocks with $\cup_i\lambda_i=\{1,...,n\}$ and $l$ is the number of blocks. The number of partitions of $n$ elements is given by Bell number $B_n$, which satisfies 
\begin{align*}
    B_{n+1}=\sum_{k=0}^n\binom{n}{k}B_k,\qquad B_0=1.
\end{align*}
The number of set partitions of $\{1,...,n\}$ according to the number of blocks is given by the \textit{Stirling numbers of the second kind} $S(n,k)$, where $k$ is the number of blocks. The definition includes also falling factorials
\begin{align*}
    \sum_{k=0}^nS(n,k)x(x-1)...(x-k+1)=x^n.
\end{align*}
Stirling numbers of the first kind and Stirling numbers of the second kind interpreted as matrix are inverse to each other 
\begin{align}
    \sum_{q}S(n,q)s(q,k)=\sum_{q}s(n,q)S(q,k)=\delta_{n,k},
\end{align}
with $S(n,k)=s(n,k)=0$ if $n<k$.
Let us refine the enumeration of set partitions a bit further as above for the permutations. Let $[a]=(a_1,a_2,...,a_k)\vdash n$ be again an integer partition of the integer $n$. Then, the number of partitions of $\{1,...,n\}$ with blocks of size $[a]\vdash n$ is
\begin{align*}
    C_{n,[a]}=\frac{n!}{\mathrm{sym}(a)\prod_ia_i!},
\end{align*}
where $\mathrm{sym}(a)$ is the same symmetry factor as before. 

Forgetting the cyclic order of the cycles of a permutation and considering these as blocks, one obtains a partition from a permutation. In other words, taking permutations modulo the cyclic order, these equivalence classes are given by partitions.

Stirling numbers of the second kind and Bell numbers are computed from $C_{n,[a]}$ via 
\begin{align*}
    S(n,k)=\sum_{\substack{[a]\vdash n\\ |[a]|=k}}C_{n,[a]},\qquad B_n=\sum_{\substack{[a]\vdash n}}C_{n,[a]}.
\end{align*}
Now let $\kappa_i$ be indeterminates, we can define the following moments
\begin{align*}
    m_n=\sum_{[a]\vdash n}C_{n,[a]}\prod_i\kappa_{a_i},
\end{align*}
which associates $\kappa_i$ to a block of length $i$ and sums over all integer partitions $[a]\vdash n$. In probability theory or statistical mechanics, the $\kappa_i$'s are called the connected parts of the moment $m_n$.

Let $ X(y)$ be again the generating series of $\kappa_i$ as above \eqref{GenKappa}, then define the generating series of partitions $W_{par}(x)$ via
\begin{align}
    W_{par}(x)=&\frac{1}{x}+\sum_{n=1}\frac{m_n}{x^{n+1}}=\frac{1}{x}+\sum_{n=1}\frac{\sum_{[a]\vdash n}C_{n,[a]}\prod_i\kappa_{a_i}}{x^{n+1}}.
\end{align}
Again, one might be interest of the relation between $ X(y)$ and $W_{par}(x)$. 

We list again three examples:
\begin{align*}
    &\text{Bell numbers:} &&\kappa_i=1,\quad  X(y)=\frac{1}{y}+\frac{1}{1-y},\quad W_{par}(x)=\sum_{n=0}^\infty \frac{B_n}{x^{n+1}}\\
    &\text{Stirling 2nd:} &&\kappa_i\kappa,\quad X(y)=\frac{1}{y}+\frac{\kappa}{1-y},\quad W_{par}(x)=\sum_{n=0}^\infty \frac{\sum_{k=0}^nS(n,k) \kappa^k}{x^{n+1}}\\
    &\text{Harer-Zagier:} &&\kappa_i=\delta_{i,2},\quad X(y)=\frac{1}{y}+y,\quad\quad W_{par}(x)=\sum_{n=0}^\infty \frac{(2n-1)!!}{ x^{2n+1}}=W_{per}(x).
\end{align*}
For the example of Harer-Zagier, the partitions and permutations coincide since we are just allowing cycles of lengths 2 for the permutations, which are in one-to-one correspondence to the blocks of the set partition. Note that the generating series of the set partitions $W_{par}(x)$ is a formal power series in $\frac{1}{x}$, which does not converge in general.

From the generating series $W_{per}(x)$ and  $W_{par}(x)$ defined above, the coefficients are computed via
\begin{align}\label{res1}
    \alpha_n=&-\Res_{x\to \infty} W_{per}(x) x^n dx\\\label{res2}
    m_n=&-\Res_{x\to \infty} W_{par}(x) x^n dx
\end{align}

\subsection{Genus expansion of permutations and partitions}
\begin{figure}[htb]
{\hspace*{2ex}
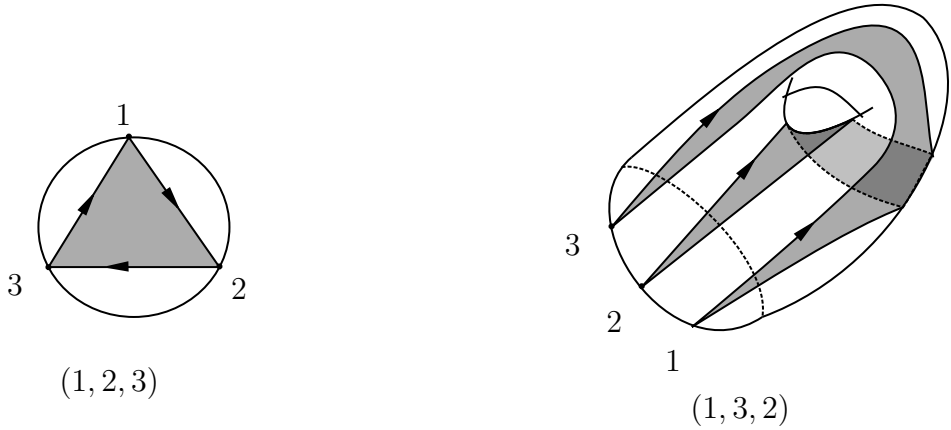}
\caption{Two permutations $\sigma_1,\sigma_2\in \mathcal{S}_3$ with $\sigma_1=(1,2,3)$ and $\sigma_2=(1,3,2)$ are drawn together with its embedding on the corresponding Riemann surface. The permutation $\sigma_1$ is of genus $g=0$ and the permutation $\sigma_2$ of genus $g=1$. The orientation indicated by the arrows matters and the enclosed area is to the right of the orientation.}
\label{Fig:k3g1}
\end{figure}

\begin{figure}[htb]
{\hspace*{2ex}
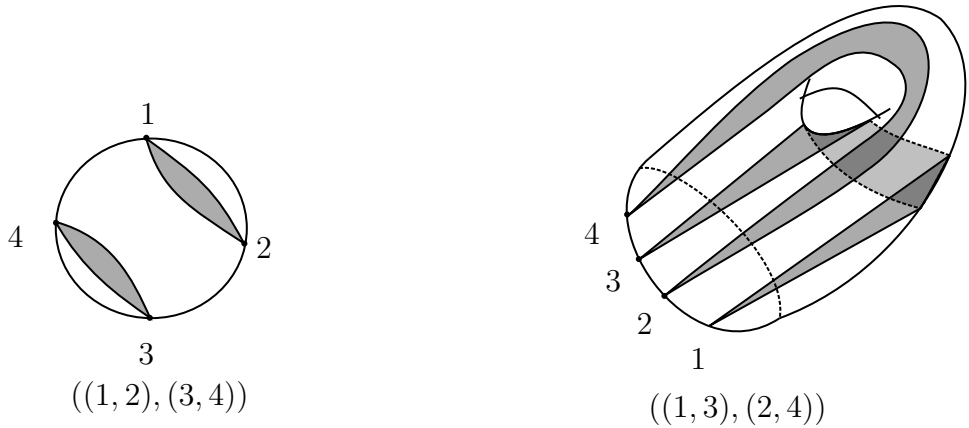}
\caption{Two partitions $\lambda,\lambda'\in P(4)$ with $\lambda=\{(1,2),(3,4)\}$ and $\lambda'=\{(1,3),(2,4)\}$ are drawn together with their embedding on the corresponding Riemann surfaces. The partition $\lambda$ is of genus $g=0$ and the partition $\lambda'$ of genus $g=1$. There is no orientation, the blocks of $\lambda,\lambda'$ are sets.}
\label{Fig:k22g1}
\end{figure}

Modify the counting problem by decomposing permutations and partitions wrt their genus was considered for example in \cite{MR3261809,MR3891093,zuber2023counting}, see also \cite{MR1649966,MR2550161} for permutations. There is a canonical way to assign a genus toa given permutation and a given partition.
This means graphically that a permutation or partition can be drawn on a Riemann surface of genus $g$ with one boundary with $n$ points on the boundary representing the set $\{1,...,n\}$ in consecutive order. A permutation connects the points on the boundary described by the cycle type respecting the orientation. The genus of a permutation is then the minimal 
topological genus of the Riemann surface such that all drawn cycles within the Riemann surface do not cross, see Fig. \ref{Fig:k3g1}. 
Similarly, a partition connects the points on the boundary described by the blocks (no orientation). The genus of a partition is the minimal topological genus of the Riemann surface such that all blocks within the Riemann surface are non-crossing, see Fig. \ref{Fig:k22g1}. 

More formally, the genus $g(\sigma)$ of a permutation $\sigma\in \mathcal{S}_n$ is defined via \cite{MR3261809}
\begin{align}\label{genus}
    2g(\sigma)=n+1-l(\sigma)-l(\sigma^{-1}\zeta_n),
\end{align}
where $l(\sigma')$ for $\sigma'\in \mathcal{S}_n$ is the number of disjoint cycles of $\sigma'$. Further, $\zeta_n\in \mathcal{S}_n$ is the circular permutation $\zeta_n(i)=i+1$ modulo $n$ and $\sigma^{-1}\zeta_n$ is the composition of $\sigma^{-1}$ with $\zeta_n$.

The genus of a partition $\lambda\in P(n)$ is defined with the same formula by considering the partition $\lambda$ as a permutation $\sigma\in \mathcal{S}_n$, where the blocks of $\lambda$ are turned into cycles of $\sigma$ with the order of the cycles constructed as increasing lists of integers. 

We expand the generating series of permutations and partitions wrt the genus
\begin{align*}
    W_{per}(x)=&\sum_{g=0}^\infty W^{(g)}_{per}(x)\\
    W_{par}(x)=&\sum_{g=0}^\infty W^{(g)}_{par}(x),
\end{align*}
where $W^{(g)}_{per}(x)$ and $W^{(g)}_{par}(x)$ are the generating series of permutations or partitions of genus $g$, respectively. Extending this notation to the coefficients, we define
\begin{align*}
    W^{(g)}_{per}(x)=&\frac{\delta_{g,0}}{x}+\sum_{n=1}\frac{\alpha^{(g)}_n}{x^{n+1}}, \qquad \alpha^{(g)}_n=\sum_{[a]\vdash n}D^{(g)}_{n,[a]}\prod_i\kappa_{a_i}\\
    W^{(g)}_{par}(x)=&\frac{\delta_{g,0}}{x}+\sum_{n=1}\frac{m^{(g)}_n}{x^{n+1}}, \qquad m^{(g)}_n=\sum_{[a]\vdash n}C^{(g)}_{n,[a]}\prod_i\kappa_{a_i},
\end{align*}
where $D^{(g)}_{n,[a]}$ counts the genus $g$ permutations of type $[a]\vdash n$ and $C^{(g)}_{n,[a]}$ counts the genus $g$ partitions of type $[a]\vdash n$.

It is very common to add a further formal parameter $\hbar$ separating the different genera
\begin{align}
    \mathcal{W}_{per}(x)=&\sum_{g=0}^\infty\hbar^{2g-1} W^{(g)}_{per}(x)\\
    \mathcal{W}_{par}(x)=&\sum_{g=0}^\infty\hbar^{2g-1} W^{(g)}_{par}(x),
\end{align}
with equality $\mathcal{W}_{per}\vert_{\hbar=1}=W_{per}$ and $\mathcal{W}_{par}\vert_{\hbar=1}=W_{par}$ at the level of formal expressions. 

For any genus and any particular choice of $\kappa_i$, the genus generating functions $W^{(g)}_{per}(x)$ and $W^{(g)}_{per}(x)$ do converge. The convergent radius depends on the $\kappa_i$ and is independent of the genus. More interestingly, evaluating the genus generating functions at $X(y)$ which was defined in \eqref{GenKappa}, one achieve a rational function in $y$, i.e. 
\begin{align*}
    W^{(g)}_{per}(X(y)),\qquad W^{(g)}_{par}(X(y))
\end{align*}
are rational in $y$. We provide for any genus $g$ the explicit formula of the generating function for permutations $W^{(g)}_{per}(X(y))$ in terms of derivatives of $X(y)$. The provided results are special cases of recent developments in the theory of Free Probability \cite{Borot:2021thu} and it relation to the theory of Topological Recursion \cite{Eynard:2007kz}. 

\begin{example}
    For the Harer-Zagier example $X(y)=\frac{1}{y}+y$, the coefficients $ \alpha_n^{(g)}=m_n^{(g)}$ correspond to the coefficients $\varepsilon_g(n)$ in \cite{MR848681}, which gives us a way to compute the virtual Euler characteristic of the moduli space of complex curves, see also \cite{MR2036721} for a comprehensible exposition. 
\end{example}

The generating series of partitions up to genus $g=2$ was derived in \cite{zuber2023counting}. We rewrite the formulas in our notation and prove some conjectures stated in \cite{coquereaux2023counting} via residue computations. More precisely, the coefficients $\alpha_n^{(g)}$ and $m_n^{(g)}$ can be computed easily from the rational functions $W^{(g)}_{per}(X(y))$ and $W^{(g)}_{par}(X(y))$ respectively via the residue formula
\begin{align}\label{res1g}
    \alpha_n^{(g)}=&-\Res_{y\to 0} W_{per}^{(g)}(X(y)) X(y)^n X'(y) dy\\\label{res2g}
    m_n^{(g)}=&-\Res_{y\to 0} W_{par}^{(g)}(X(y)) X(y)^n X'(y) dy,
\end{align}
where $X(y)$ was defined in \eqref{GenKappa}.

\begin{remark}
    The notion of genus permutations used in this article is equivalent to so-called unicellular maps which are well studied in the combinatorics literature, see for instance \cite{CHAPUY2011874,CHAPUY20132064}. 
\end{remark}

\begin{remark}\label{rem1}
    The background of the this article concerning the recent developments in the theory of Free Probability and its connection to the $x-y$ symplectic transformation within the theory of Topological Recursion is beyond the scope of this article. We refer an ambitious reader to the growing literature \cite{Bychkov:2020ujd,Bychkov:2020yzy,Bychkov:2021hfh,Alexandrov:2022ydc,Borot:2021eif,Borot:2021thu, Hock:2022wer,Hock:2022pbw} and references therein. The explicit formula for $W^{(g)}_{per}(X(y))$ of Sec. \ref{sec2} is a prime example which did not appear (to our knowledge) in the literature, yet. It is a special case of a much more general concept, but there is no need to recap the general formulae and all of its technicalities.  
\end{remark}

\subsection{Planar Permutations and Partitions}\label{Sec:g0}
The genus $g=0$ permutations or partitions which are also called non-crossing permutations or non-crossing partitions appeared first in \cite{MR309747}. They reappeared in different contexts coming from physics \cite{Brezin:1977sv,Cvitanovic:1980jz} or operator algebras \cite{Voiculescu1986AdditionOC}. The work of Voiculescu gave rise to the theory of Free Probability. 

At genus $g=0$, permutations and partitions coincide $W^{(0)}_{per}(x)=W^{(0)}_{par}(x)$, this is a well-known fact and reproved in the literature at several places  \cite{MR0404045,10.1155/S1073792804133023}. Using Langrange inversion formula (see for instance \cite{Brezin:1977sv}), it is easy to prove that $W^{(0)}_{per}(x)$ is the formal inverse of $X(y)$, i.e.
\begin{align*}
    W^{(0)}_{per}(X(y))=y,\qquad X( W^{(0)}_{per}(x))=x.
\end{align*}
The special case of $\kappa_i=1$, i.e. $X(y)=\frac{1}{y(1-y)}$ is counted by the Catalan numbers.

\section{Genus Permutations}\label{sec2}
\subsection{Closed Formula for the Generating Function of Permutations at all Genera}
Genus permutations were consider recently in a much broader sense than described here in this article \cite{Borot:2021thu}. The motivation came from Free Probability \cite{Voiculescu1986AdditionOC,Collins2006SecondOF} to give a general formula between the so-called higher order free cumulants and higher order moments, see also Remark \ref{rem1}.  This article explains explicit formulae for a special case considered in \cite{Borot:2021thu}. We will not include higher order free cumulants, i.e. we put all $\kappa_{i_1,i_2,...,i_n}=0$ for $n\geq 2$ appearing in the article \textit{op. cit}. These $\kappa_{i_1,i_2,...,i_n}$ carry by themselves topological structure, which means that our example is reduced to disc-like cycles for permutations, which are exactly the ones described above.

For this setting, we have the following explicit expression for all $W^{(g)}_{per}(X(y))$:

\begin{theorem}\label{ThmPerm}
    For any $X(y)=\frac{1}{y}+\sum_{i\geq 1}\kappa_i y^{i-1}$, the generating function of permutations of genus $g$, where $i$-cycles are weighted by $\kappa_i$, is given by
\begin{align*}
   W^{(g)}_{per}(X(y))=\sum_{l=1}^g\sum_{\substack{[a]\vdash g\\ [a]=(a_1,...,a_l)}}\bigg(-\frac{d}{d X(y)}\bigg)^{2g+l-1}\frac{\bigg(-\frac{1}{X'(y)}\prod_{i=1}^l\frac{X^{(2 a_i)}(y)}{2^{2 a_i}(2 a_i+1)!}\bigg)}{\mathrm{sym}(a)},
\end{align*}
where $\mathrm{sym}(a)$ is the symmetry factor $\prod_{j}k_j!$ with $k_1,k_2,..$ the numbers of equal parts of the integer partition $[a]$. 

Equivalently, as a formal expansion in $\hbar$, the following formula holds for $\mathcal{W}_{per}(x)=\sum_{g=0}^\infty\hbar^{2g-1} W^{(g)}_{per}(x)$:
\begin{align}\label{Whsum}
    \mathcal{W}_{per}(X(y))=\sum_{m\geq 0} \bigg(-\frac{\partial}{\partial {X(y)}}\bigg)^m\bigg(-\frac{dy}{dX(y)}\bigg)
	[u^m]\frac{\exp\bigg(\frac{\int_{y-\hbar u/2
    }^{y+\hbar u/2} X(y') dy'}{\hbar}-X(y)u\bigg)}{\hbar u}.
\end{align}
\begin{proof}
The definition of the genus \eqref{genus} is a special case of the definition used in \cite[Def. 4.4]{Borot:2021thu} restricted to one boundary. Therefore, the functional relation \cite[Thm. 3.4]{Borot:2021thu} applies. This functional relation was simplified in \cite{Hock:2022pbw} and yields for $n=1$ and unramified $y$ (which coincides with the considered problem in this article)
    \begin{align*}
    W^{(g)}_{per}(X(y))=[\hbar^{2g-1}]\sum_{m\geq 0} \bigg(-\frac{\partial}{\partial {X(y)}}\bigg)^m\bigg(-\frac{dy}{dX(y)}\bigg)
	[u^m]\frac{\exp\bigg(\hat{\Phi}^\vee_1(y;\hbar,u)-X(y) u\bigg)}{\hbar u},
\end{align*}
where 
\begin{align*}
    \hat{\Phi}^\vee_1(y;\hbar,u):=&\frac{1}{\hbar}\bigg(\Phi(y+\frac{\hbar u}{2})-\Phi(y-\frac{\hbar u}{2})\bigg),\\
    \Phi(y)=&\int_{o
    }^y X(y') dy'.
\end{align*}
This gives us the formula for $\mathcal{W}_{per}(x)$. 

The explicit structure of the coefficients in this formula at each order in $\hbar^{2g-1}$ was analysed in \cite[\S 2.4]{Hock:2022pbw}, where we have to consider just the special case of one $\circ$-vertex (corresponding to $n=1$) and just 1-valent $\bullet$-vertices with weight 0 (due to the fact that $y$ is unramified in the sense of \cite{Hock:2022pbw}), but arbitrary edge weights. The coefficients are literally the coefficients of $S(t)=\frac{e^{t/2}-e^{-t/2}}{t}=\sum_{k=0}^\infty \frac{t^{2k}}{2^{2k}(2k+1)!}=1+\frac{t^2}{24}+\frac{t^4}{1920}+\frac{t^6}{322560}+...$ 

Expanding the exponential, the coefficients of the formal expansion in $u$ and $\hbar$ are not independent. The leading order of the argument of the exponential is $\hbar$. Collecting the same orders of the expanded exponential at the order $\hbar^{2g-1}$ has different contributions from different orders of the argument of the exponential.  
It breaks down to a sum of integer partitions of $g$, i.e. $[a]\vdash g$ with $[a]=(a_1,...,a_l)$ where $a_i\in \mathbb{N}_{>0} $ with $\sum_{i=1}^l a_i=g$ and $l$ the length of the integer partition $[a]$.
\end{proof}
\end{theorem}
This explicit functional relation for genus permutations did not appear in literature up to our knowledge. There exists some work on genus permutations \cite{MR1649966,MR2550161,CHAPUY20132064} which does not consider generating series and neither functional relations. The functional relation of Theorem \ref{ThmPerm} is therefore more in the spirit of the original works \cite{MR309747,Brezin:1977sv,Cvitanovic:1980jz}.

\begin{example}\label{experg}
The examples for $g\leq 3$ are listed:
    \begin{align*}
    W^{(0)}_{per}(X(y))=&y,\\
     W^{(1)}_{per}(X(y))=&-\frac{d^2}{d X(y)^2}\bigg(\frac{1}{X'(y)}\frac{X''(y)}{24}\bigg),\\
      W^{(2)}_{per}(X(y))=&\frac{d^5}{d X(y)^5}\bigg(\frac{1}{X'(y)}\frac{(X''(y))^2}{24^2\cdot 2}\bigg)-\frac{d^4}{d X(y)^4}\bigg(\frac{1}{X'(y)}\frac{X^{(4)}(y)}{1920}\bigg),\\
      W^{(3)}_{per}(X(y))=&-\frac{d^8}{d X(y)^8}\bigg(\frac{1}{X'(y)}\frac{(X''(y))^3}{24^3\cdot 3!}\bigg)+\frac{d^7}{d X(y)^7}\bigg(\frac{1}{X'(y)}\frac{X^{(4)}(y) X^{(2)}(y)}{1920\cdot 24}\bigg)\\
      &-\frac{d^6}{d X(y)^6}\bigg(\frac{1}{X'(y)}\frac{X^{(6)}(y) }{322560}\bigg).
\end{align*}
\end{example}

\begin{example}
Taking the Example \ref{experg}, it is very easy to compute with the formula \eqref{res1g} 
\begin{align*}
    \alpha_n^{(g)}=-\Res_{y\to 0} W^{(g)}_{per}(X(y)) X(y)^n X'(y) dy
\end{align*}
the coefficients $\alpha_n^{(g)}$ for general $\kappa_i$'s. We get for instance: 
    \begin{align*}
        \alpha_{3}^{(1)}=&\kappa_3\\
        \alpha_{4}^{(1)}=&4 \kappa_{3} \kappa_{1}+\kappa_{2}^{2}+5 \kappa_{4}\\
         \alpha_{5}^{(1)}=&10 \kappa_{1}^{2} \kappa_{3}+5 \kappa_{1} \kappa_{2}^{2}+25 \kappa_{1} \kappa_{4}+15 \kappa_{2} \kappa_{3}+15 \kappa_{5}\\
          \alpha_{6}^{(1)}=&20 \kappa_{1}^{3} \kappa_{3}+15 \kappa_{1}^{2} \kappa_{2}^{2}+75 \kappa_{1}^{2} \kappa_{4}+90 \kappa_{1} \kappa_{2} \kappa_{3}+10 \kappa_{2}^{3}+90 \kappa_{1} \kappa_{5}+60 \kappa_{2} \kappa_{4}+25 \kappa_{3}^{2}+35 \kappa_{6}\\
          \alpha_{5}^{(2)}=&8 \kappa_{5}\\
        \alpha_{6}^{(2)}=&48 \kappa_{5} \kappa_{1}+24 \kappa_{2} \kappa_{4}+12 \kappa_{3}^{2}+84 \kappa_{6}\\
         \alpha_{7}^{(2)}=&168 \kappa_{1}^{2} \kappa_{5}+168 \kappa_{1} \kappa_{2} \kappa_{4}+84 \kappa_{1} \kappa_{3}^{2}+49 \kappa_{2}^{2} \kappa_{3}+588 \kappa_{1} \kappa_{6}+322 \kappa_{2} \kappa_{5}+273 \kappa_{3} \kappa_{4}+469 \kappa_{7}\\
         \alpha_{7}^{(3)}=&180 \kappa_{7}\\
         \alpha_{8}^{(3)}=&1440 \kappa_{1} \kappa_{7}+720 \kappa_{2} \kappa_{6}+608 \kappa_{3} \kappa_{5}+276 \kappa_{4}^{2}+3044 \kappa_{8}
    \end{align*}
\end{example}
For any $g>0$, the function $W^{(g)}_{per}(X(y))$ is a rational function in $y$ with poles just located at the ramification points of $X(y)$, i.e. $\{y\in \mathbb{C}| X'(y)=0\}$. This can be proved easily from the Theorem. The leading pole for $W^{(g)}_{per}(X(y))$ is of order $6g-1$, which is generated if $[a] \vdash g$ is $[a]=(\underbrace{1,...,1}_{g})$ and all derivatives $\frac{\partial}{\partial X(y)}$ acting on $\frac{1}{X'(y)}$. As an expansion about $x=\infty$, $W^{(g)}_{per}(x)\in \mathcal{O}(\frac{1}{x^{2g+1}})$. 

Note that $W^{(g)}_{per}(x)$ makes sense as a formal expansion at $x=\infty$, whereas $W^{(g)}_{per}(X(y))$ is an analytic continuation defined for $y$ in the complex plane. Constructing the inverse of $X(y)$ is just locally possible. This prevents us to write down explicit formulas for $W^{(g)}_{per}(x)$ in general.  

If $X(y)$ is a ramified covering of degree two, the inverse can be written explicitly. Choosing the correct branch it is possible to write down explicit formulae for $W^{(g)}_{per}(x)$. This type of examples are considered in Appendix \ref{appper} by specifying $\kappa_i$.

\subsection{Laplace transform}
The formal Laplace transform is a very natural operation on the formula \eqref{Whsum} of Theorem \ref{ThmPerm}. In a more general setting, the Laplace transform was already considered in \cite{Hock:2023qii} deriving new formulas for intersection numbers on the moduli space of complex curves. We will recap the few computational steps due to the important insight they are bringing. We are not precise about the integration contour which should be chosen in the complex plane such that computation converges (see \textit{op. cit.} for details):
\begin{theorem}[\cite{Hock:2023qii}]
    The formal Laplace transform of $\mathcal{W}_{per}(x)=\sum_{g=0}^\infty\hbar^{2g-1} W^{(g)}_{per}(x)$ is given by
    \begin{align}\label{Laplace}
        \int dx e^{-x \mu }\mathcal{W}_{per}(x)=\int \frac{dy}{\hbar \mu}\bigg(\frac{y-\frac{\hbar\mu}{2}}{y+\frac{\hbar\mu}{2}}\bigg)^{1/\hbar} \exp\bigg(\frac{\sum_{i\geq 1}\frac{\kappa_i}{i}[(y-\frac{\hbar\mu}{2})^i-(y+\frac{\hbar\mu}{2})^i]}{\hbar}\bigg).
    \end{align}
    \begin{proof}
    Assume the integration contour vanishes all boundary terms and as a formal expansion in $\hbar$ converges for any coefficient. We compute directly the formal Laplace transform of \eqref{Whsum} with $x=X(y)$
        \begin{align*}
            &\int dX(y) e^{-X(y) \mu }\mathcal{W}_{per}(X(y))\\
            =&\int dX(y) e^{-X(y) \mu }\sum_{m\geq 0} \bigg(-\frac{\partial}{\partial {X(y)}}\bigg)^m\bigg(-\frac{dy}{dX(y)}\bigg)
	[u^m]\frac{\exp{\bigg(\frac{\int_{y-\hbar u/2
    }^{y+\hbar u/2} X(y') dy'}{\hbar}-X(y)u\bigg)}}{\hbar u}\\
    =&\int dX(y) e^{-X(y) \mu }\sum_{m\geq 0} (-\mu)^m\bigg(-\frac{dy}{dX(y)}\bigg)
	[u^m]\frac{\exp{\bigg(\frac{\int_{y-\hbar u/2
    }^{y+\hbar u/2} X(y') dy'}{\hbar}-X(y)u\bigg)}}{\hbar u}\\
    =&\int dX(y) e^{-X(y) \mu }\bigg(-\frac{dy}{dX(y)}\bigg)
	\frac{\exp{\bigg(\frac{\int_{y-\hbar u/2
    }^{y-\hbar \mu/2} X(y') dy'}{\hbar}+X(y)\mu\bigg)}}{-\hbar \mu}\\
    =&\int \frac{dy}{\hbar \mu} \exp{\bigg(\frac{\int_{y-\hbar u/2}^{y-\hbar \mu/2} X(y')dy'}{\hbar}\bigg)},
        \end{align*}
        which is equivalent to the claimed result. The following computational steps are performed
        \begin{itemize}
            \item Inserting \eqref{Whsum}
            \item Integrating by parts $m$ times for each summand in the $m$-summation
            \item Using for a formal expression $P(\mu)=\sum_m \mu^m [u^m]P(u)$
            \item Change of variables and cancellation of the Laplace kernel
            \item Integrating $X(y)=\frac{1}{y}+\sum_{i\geq 1}\kappa_i y^{i-1}$ in the exponential.
        \end{itemize}
    \end{proof}
\end{theorem}
The important insight is that the formal Laplace transform seems to be the canonical object to look at. For the rhs of \eqref{Laplace}, the integrand is understood to be formally expanded in $\hbar$ before carrying out the integration.

\subsection{Generalisations}\label{sec:General}
\begin{figure}[htb]
{\hspace*{2ex}
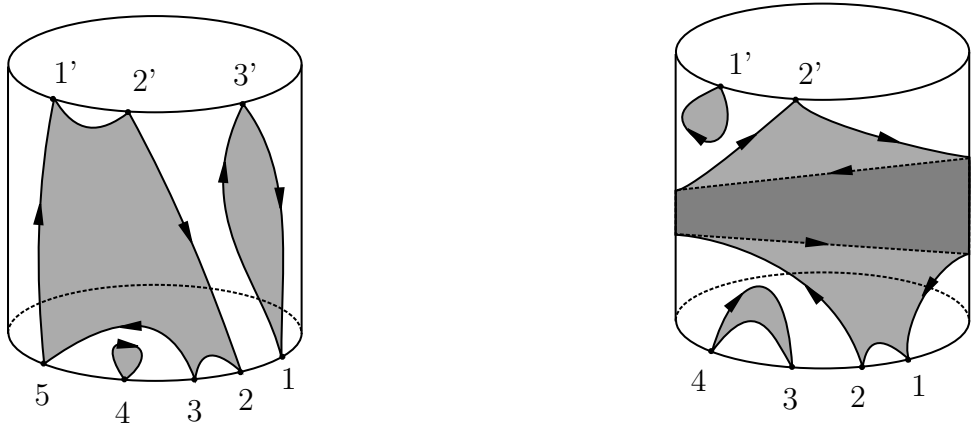}
\caption{On the left, we have the boundary permutation $\tau\in \mathcal{S}_{8}$ and the permutation  $\sigma\in \mathcal{S}_{8}$ of genus $g=0$ of the form $\tau=((1,2,3,4,5),(1',2',3'))$ and  $\sigma=((1,3'),(2,3,5,1',2'),(4))$. On the right, we have a partitioned permutation, where the associated weight is of the form $\kappa_1\kappa_2\kappa_{1,2}$, where $\kappa_1$ comes from the cycle $(1')$, $\kappa_2$ from $(3,4)$ and $\kappa_{1,2}$ from $(2'),(1,2)$ connected non-trivially.}
\label{Fig:cylinder}
\end{figure}

There are two types of generalisations which can be made. Both are included in the more general theory of partitioned permutations \cite{Collins2006SecondOF} and/or surfaced permutations considered and understood in general in \cite{Borot:2021thu}, see also \cite{Kazarian:2014ooa}. We distinguish them by
\begin{itemize}
    \item Allowing the Riemann surface to have more than one boundary
    \item Allowing permutations to carry additional topological structure.
\end{itemize}
If we have $b$ boundaries, we associate to the set $\{1,...,n\}$ a permutation $\tau\in \mathcal{S}_n$ encoding the structure of the boundary in the sense that $\tau$ has $b$ cycles, each associated to a boundary component.  Then, there is a second permutation $\sigma\in\mathcal{S}_n$ on the same set such that the boundaries encoded by $\tau$ are connected through $\sigma$, see Fig.\ref{Fig:cylinder}. One can associate a genus to the pair $(\tau,\sigma)$ via (see for instance \cite{MR2036721,MR3261809})
\begin{align}\label{genusb}
    2g(\sigma,\tau)=n+2-b-l(\sigma)-l(\sigma^{-1} \tau ),
\end{align}
where $l(\sigma),b=l(\tau)$ are the number of disjoint cycles of $\sigma,\tau\in \mathcal{S}_n$, respectively. 

The second generalisation includes partitioned permutations in the sense of \cite{Collins2006SecondOF}. Permutations themselves can carry a topological structure. The right figure of Fig. \ref{Fig:cylinder} shows an example of a partitioned permutation having additional topological structure of a cylinder connecting $(2',1,2)$ topologically non-trivially.  
We will just review the cylinder case (understood in \textit{op. cit.} and we refer to this reference for further details). Associate a weight of the form $\kappa_{i,j}$ to a topologically non-trivial cycle related to a second order free cumulant. Let the generating series of $\kappa_{i,j}$ be
\begin{align*}
    X(y_1,y_2)=\sum_{i,j\geq 1}\kappa_{i,j}y_1^{i-1}y_2^{j-1},
\end{align*}
then the generating function of $W^{(0)}_{perm,2}(x_1,x_2)=\sum_{i,j\geq 1}\frac{\alpha^{(0)}_{i,j}}{x_1^{i+1}x_2^{j+1}}$ of permutations on the cylinder including second order free cumulants is given via \cite[Theorem 6.3]{Collins2006SecondOF}\footnote{the relation between the different generating series between this article and \cite{Collins2006SecondOF} is as follows: $C(x)\to \frac{X(y)}{y},M(x)\to  \frac{W^{(0)}_{perm}(\frac{1}{x})}{x}, C(x,y)\to \frac{X(y_1,y_2)}{y_1y_2}$ and $M(x,y)\to  \frac{W^{(0)}_{perm,2}(\frac{1}{x_1},\frac{1}{x_2})}{x_1x_2}$}
\begin{align}\label{W02Perm}
    W^{(0)}_{perm,2}(X(y_1),X(y_2))=\frac{X(y_1,y_2)+\frac{1}{(y_1-y_2)^2}}{X'(y_1)X'(y_2)}-\frac{1}{(X(y_1)-X(y_2))^2}.
\end{align}
We mention this result here to use it later to derive a new result for the cylinder moments of non-crossing partitions in Sec. \ref{Sec.PartBound}. The first explicit coefficients $\alpha^{(0)}_{i,j}$ are listed in \cite[Remark 6.6]{Collins2006SecondOF}\footnote{There is a typo for the $\alpha_{1,3}$: the term $3\kappa_1^2\kappa_{1,1}$ is missing}.

\begin{remark}\label{rem:gbperm}
    The important result of \cite{Borot:2021thu} is a tremendous generalisation including higher genus and higher order cumulants of the form $\kappa_{n_1,...,n_b}^{(g)}$ and relate them to the moments
\begin{align*}
    \alpha^{(g)}_{n_1,...,n_b}.
\end{align*}
Collecting these in a generating series of the form
\begin{align*}
    W^{(g)}_{per,b}(x_1,...,x_b)=\frac{\delta_{g,0}\delta_{b,1}\delta_{n_1,1}}{x}+\sum_{n_1,...,n_b\geq 1}\frac{\alpha^{(g)}_{n_1,...,n_b}}{x_1^{n_1+1}...x_b^{n_b+1}},
\end{align*}
the explicit relation to the generating series of the $\kappa_{n_1,...,n_b}^{(g)}$ is known, which is formulated in \cite[Theorem 1.1]{Borot:2021thu} and simplified in \cite{Hock:2022pbw}, but it is beyond the scope of the present article. Furthermore, the formulae presented in the Thm \ref{ThmPerm} are even further simplifications and specialisations.
\end{remark}

\begin{remark}
    This remark makes the connection to the theory of Topological Recursion \cite{Eynard:2007kz}. Taking the spectral curve as defined in \textup{op.cit.} to be $(\mathbb{P}^1,x=\frac{1}{z}+\sum_{i\geq 1}\kappa_{i}z^{i-1},y=z,\frac{dz_1 dz_2}{(z_1-z_2)^2})$, then all $\omega_{g,n}$ generated by the algorithm of Topological Recursion are essentially $W^{(g)}_{per,b}$, where all cumulants $\kappa_{n_1,...,n_b}^{(g)}=0$ except the $\kappa_i$'s. Using the formula of Topological Recursion to compute $W^{(g)}_{per}(x)$ is almost impossible since the branched covering $x$ has arbitrary high degree. However, the recent insight \cite{Borot:2021thu,Hock:2022pbw,Alexandrov:2022ydc} of $x-y$ symplectic transformation for Topological Recursion provides the compact result of Theorem \ref{ThmPerm} and gives furthermore a formula for all other $W^{(g)}_{per,b}$ with $b>1$. 
\end{remark}

\section{Genus Partitions}
Counting genus partitions wrt the sizes of the blocks turns out to be completely different. The idea of forgetting the cyclic structure of a permutation to get a partition has topologically nontrivial consequences. For genus $g=0$, the number of permutations and partitions is equal (see Sec. \ref{Sec:g0}) since every cycle has a unique ordering by restricting to $g=0$. This changes drastically at higher genus and including more boundaries.  Genus partitions were considered for instance in \cite{MR3261809,MR3891093,zuber2023counting,coquereaux2023counting}.

\subsection{Closed Formula for the Generating Function of Partitions for $g<3$}
This subsection brings per se no new results on genus partitions. We will rewrite the result of \cite{zuber2023counting} in our notation, which seems to be the canonical representation. The benefit of this representation is that we can prove easily (see Appendix \ref{app:part}) some conjectures made about specific genus partitions in \cite{coquereaux2023counting}. We will use the following representation for this:
\begin{proposition}\label{PropPart}
    Let $X(y)=\frac{1}{y}+\sum_{i\geq 1}\kappa_i y^{i-1}$, then the number of genus partitions according to the sizes of their blocks are counted for genus 1 by
    \begin{align*}
        W^{(1)}_{par}(X(y))=&\frac{\partial}{\partial X(y)}\bigg(\frac{1}{4 y^4X'(y)^2}+\frac{1}{6 y^6X'(y)^3}\bigg)
        \end{align*}
        and for genus 2 by
        \begin{align*}
        W^{(2)}_{par}(X(y))=&\frac{\partial}{\partial X(y)}\bigg[ \frac{21}{8 y^{8} X'(y)^{4}}+\frac{74}{5 y^{10} X'(y)^{5}}+\frac{24}{y^{12} X'(y)^{6}}+\frac{12}{y^{14} X'(y)^{7}}\\\nonumber
    & -\frac{X^{(3)}(y)}{8 y^{8} X'(y)^{6}} -\frac{X^{(3)}(y)}{4 y^{10} X'(y)^{7}}-\frac{X^{(3)}(y)}{8 y^{12} X'(y)^{8}}  \\\nonumber
    &  +\frac{\left(X^{(2)}(y)\right)^{2}}{24 y^{6} X'(y)^{6}}+\frac{\left(X^{(2)}(y)\right)^{2}}{y^{8} X'(y)^{7}}+\frac{19 \left(X^{(2)}(y)\right)^{2}}{8 y^{10} X'(y)^{8}}+\frac{35 \left(X^{(2)}(y)\right)^{2}}{24 y^{12} X'(y)^{9}}  \\\nonumber
    &   +\frac{X^{(2)}(y)}{y^{7} X'(y)^{5}}+\frac{23 X^{(2)}(y)}{3 y^{9} X'(y)^{6}}+\frac{29 X^{(2)}(y)}{2 y^{11} X'(y)^{7}}+\frac{8 X^{(2)}(y)}{y^{13} X'(y)^{8}} \bigg].
    \end{align*}
    \begin{proof}
        This is just a rewriting of the results of \cite{zuber2023counting} by identifying the generating series
\begin{align*}
    W^{(g)}_{par}(x)=&\frac{Z^{(g)}(\frac{1}{x})}{x},\qquad 
    X(y)=\frac{1+W(y)}{y}.
\end{align*}

Applying this definitions to the $\tilde{X}_l(x)$,  $\tilde{Y}_l(x)$ and $V(x)$ of \cite{zuber2023counting} gives us the following expressions in terms of $X(y)$:
\begin{align*}
    \tilde{X}_l(x)&\quad \rightarrow\quad \frac{\frac{(-y)^l}{(l-1)!} X^{(l-1)}(y)+1}{y^{2l}(X'(y))^l}, \qquad l>2,\\
    \tilde{X}_2(x)&\quad \rightarrow \quad-\frac{y^2 X'(y)+1}{y^{2}X'(y)},\\
    \tilde{Y}_l(x)&\quad \rightarrow \quad \frac{(-y)^l (X(y) y)^{(l)}}{l!y^{2l}(X'(y))^l},\\
    1-V(x)&\quad \rightarrow \quad-\frac{y X'(y)}{X(y)}.
\end{align*}
Note that the denominator of $1-V(x)$ cancels exactly the additional factor of $x$ in the definition of $W^{(g)}_{par}(x)$ comparing to $Z^{(g)}(x)$. 

With all definitions in place, we write the functional relations \cite[Thm 1 \& 2]{zuber2023counting} of $Z^{(1)}(x)$ and $Z^{(2)}(x)$ in our notation. We find that a derivative wrt $X(y)$ can be pulled out in general and get the assertion.
    \end{proof}
\end{proposition}
The new representation shows that $W^{(g=1,2)}_{par}(X(y))$ has just poles at the ramification points of $X(y)$, i.e. $\{y\in \mathbb{C}| X'(y)=0\}$, the possible pole at $y=0$ vanishes due to the second order pole of $X'(y)$ at $y=0$. Furthermore, $W^{(g=1,2)}_{par}(X(y))$ has a global primitive on the complex continued $y$-plane. Thus, there is no purely first order pole at the ramification points. In other words, the residue at the ramification points vanishes. The leading pole is of order $6g-1$ for $g=1,2$. All these analyticity properties are the same for complex continued generating functions of genus permutations.

From the representation of Proposition \ref{PropPart}, the moments $m^{(g)}_{n}$ can be extracted easily:
\begin{example}
Insert $W^{(g)}_{par}(X(y))$ for generic $\kappa_i$ into the formula \eqref{res2g} 
\begin{align*}
    m_n^{(g)}=-\Res_{y\to 0} W^{(g)}_{par}(X(y)) X(y)^n X'(y) dy.
\end{align*}
First moments are given by: 
    \begin{align*}
        m_{4}^{(1)}=&\kappa _2^2\\
         m_{5}^{(1)}=&5 \kappa _1 \kappa _2^2+5 \kappa _3 \kappa _2\\
          m_{6}^{(1)}=&10 \kappa _2^3+15 \kappa _1^2 \kappa _2^2+30 \kappa _1 \kappa _3 \kappa _2+9 \kappa _4 \kappa _2+6 \kappa _3^2\\
          m_{7}^{(1)}=&35 \kappa _2^2 \kappa _1^3+105 \kappa _2 \kappa _3 \kappa _1^2+70 \kappa _2^3 \kappa _1+42 \kappa _3^2 \kappa _1+63 \kappa _2 \kappa _4 \kappa _1+70 \kappa _2^2 \kappa _3+21 \kappa _3 \kappa _4+14 \kappa _2 \kappa _5\\
          m_{6}^{(2)}=&\kappa _3^2\\
        m_{7}^{(2)}=&14 \kappa _3 \kappa _2^2+7 \kappa _1 \kappa _3^2+7 \kappa _3 \kappa _4\\
         m_{8}^{(2)}=&21 \kappa _2^4+112 \kappa _1 \kappa _3 \kappa _2^2+54 \kappa _4 \kappa _2^2+100 \kappa _3^2 \kappa _2+28 \kappa _1^2 \kappa _3^2+12 \kappa _4^2+56 \kappa _1 \kappa _3 \kappa _4+16 \kappa _3 \kappa _5
    \end{align*}
\end{example}
For $g=1,2$, the generating series are as an expansion about $x=\infty$, $W^{(g)}_{par}(x)\in \mathcal{O}(\frac{1}{x^{2g+2}})$. We prove some minor conjectures made in \cite{coquereaux2023counting} by using  Proposition \ref{PropPart} in Appendix \ref{app:part}.

\subsection{Generalisation}\label{Sec.PartBound}
For genus partitions, we can ask for the same generalisations as for permutations in Sec. \ref{sec:General}. Allowing more boundaries is already an interesting task. While for the non-crossing permutations the cylinder was understood in \cite{Collins2006SecondOF}, there is just some qualitative understanding for non-crossing partitions on the cylinder \cite{10.1155/S1073792804133023}. In the context of Coexter groups, there is some study on set partitions on a cylinder \cite{brestensky2024noncrossing} where the explicit connection to this article is not completely clear to us. We will provide the explicit formula for counting non-crossing partitions on the cylinder. On the other hand, for partitions one could also include topologically nontrivial blocks in exactly the same manner as for permutations. This means that a block of a partition can be separated in two or more parts which are connected topologically nontrivial or can even have higher genus. The topology of a given genus partition includes then all topologies of its blocks. However, we will not stress this problem in general in this article and not give precise definitions either. We want rather to provide the result for non-crossing partitions on the cylinder and postpone the general question for future work.

For the cylinder, as mentioned above, a partition can be extended with a block of cylinder topology exactly in the same way as in \cite{Collins2006SecondOF}. Forgetting the cyclic structure is trivial on the cylinder for the cylinder block. We associate a indeterminate $\kappa_{i,j}$ to a block of cylinder topology connecting $i$ points on one boundary with $j$ points on the other boundary. Define the generating function to be (the same as for permutations)
\begin{align*}
    X(y_1,y_2)=\sum_{i,j\geq 1}\kappa_{i,j}y_1^{i-1}y_2^{j-1}.
\end{align*}
Then, we find the following result for counting partitions of cylinder topology:
\begin{theorem}\label{Thm:W02Part}
Let $W^{(0)}_{par,2}(x_1,x_2)=\sum_{i,j\geq 1}\frac{m^{(0)}_{i,j}}{x_1^{i+1}x_2^{j+1}}$ be the generating series of connected partitions of cylinder type with $m^{(0)}_{i,j}$ the coefficient having $i$ elements on one boundary and $j$ on the other. A factor $\kappa_i$ is associated to a block with $i$ elements and $\kappa_{i_1,i_2}$ to a block of cylinder type connecting $i_1$ elements of one boundary with $i_2$ of the other. Let further be $X(y)=\frac{1}{y}+\sum_{i\geq 1}\kappa_iy^{i-1}$. Then, the following holds
\begin{align}\label{W02part}
W^{(0)}_{par,2}&(X(y_1),X(y_2))
=\frac{X(y_1,y_2)+\frac{1}{(y_1-y_2)^2}}{X'(y_1)X'(y_2)}-\frac{1}{(X(y_1)-X(y_2))^2}\\\nonumber
&+\frac{\frac{1}{y_1y_2}(1-y_1\partial_{y_1}y_2\partial_{y_2}) \bigg(y_1y_2\frac{X(y_1)-X(y_2)}{y_1-y_2}+\frac{1}{y_1y_2}\bigg)}{X'(y_1)X'(y_2)}.
\end{align}
\begin{proof}
The proof is a consequence of the counting problem of the permutation \eqref{W02Perm} of \cite{Collins2006SecondOF} by forgetting the cyclic order to get partitions. We have to distinguish between three contributions:
\begin{itemize}
    \item[(1)] terms including $\kappa_{i,j}$
    \item[(2)] terms including just $\kappa_{i}$'s, where more than one block connects the two boundaries
    \item[(3)] terms including just $\kappa_{i}$'s, where exactly one block connects the two boundaries.
\end{itemize}
Taking permutations into account and forgetting the cyclic structure of the cycles, the first two cases (1) and (2) are in one-to-one correspondence to partitions. 

Looking more precisely at (1), terms including cycles of cylinder type are of the form
\begin{align*}
    \kappa_{i,j}\kappa_{i_1}\kappa_{i_2}...
\end{align*}
Forgetting cyclic structure of $\kappa_{i,j}$ gives a unique partition of cylinder type, because the cycles of $\kappa_{i,j}$ are in consecutive order since otherwise it would have a higher genus due to the definition \eqref{genusb}. Furthermore, all cycles associated with $\kappa_{i_1},\kappa_{i_2},...$ give unique blocks for the partition since each cycle is just connected to one of the boundaries and we can use the same argument as for planar permutations vs partitions (disc). That is, planar permutations and partitions are in one-to-one correspondence. In Fig. \ref{Fig:cylinderPart}, the upper left permutation gives just one partition which is the lower left one.

Now we are looking at (2), which are cycles just of disc topology $\kappa_i$'s and at least two cycles are connecting the boundaries. Taking any set associated to cycles of a permutation of this type, there is no other way of generating cycles with these sets. A cycle connecting the two boundaries would need a nontrivial twist in the cylinder which would cross an other cycle connecting the two boundaries, which is not allowed. 

The last and most important case is (3), where we have just one cycle of disc topology connecting the two boundaries and all the other cycles are connected to just one boundary. This is highly related to (1), which is shown in Fig. \ref{Fig:cylinderPart}. Taking a permutation of case (1) of the form $\kappa_{i,j}\kappa_{i_1}\kappa_{i_2}...$, there are $i\cdot j$ associated permutations of the form $\kappa_{i+j}\kappa_{i_1}\kappa_{i_2}...$, where all $\kappa_{i_l}$ are the same cycles. These $i\cdot j$ possibilities arise by fixing how one element on the first boundary ($i$ possibilities) is send to an element on the second ($j$ possibilities). However, for partitions, on the other hand, there is just one way of going from $\kappa_{i,j}\to \kappa_{i+j}$, since all the $i\cdot j$ cycles of the form $\kappa_{i+j}$ give the same block in a partition. 
\begin{figure}[htb]
{\hspace*{2ex}
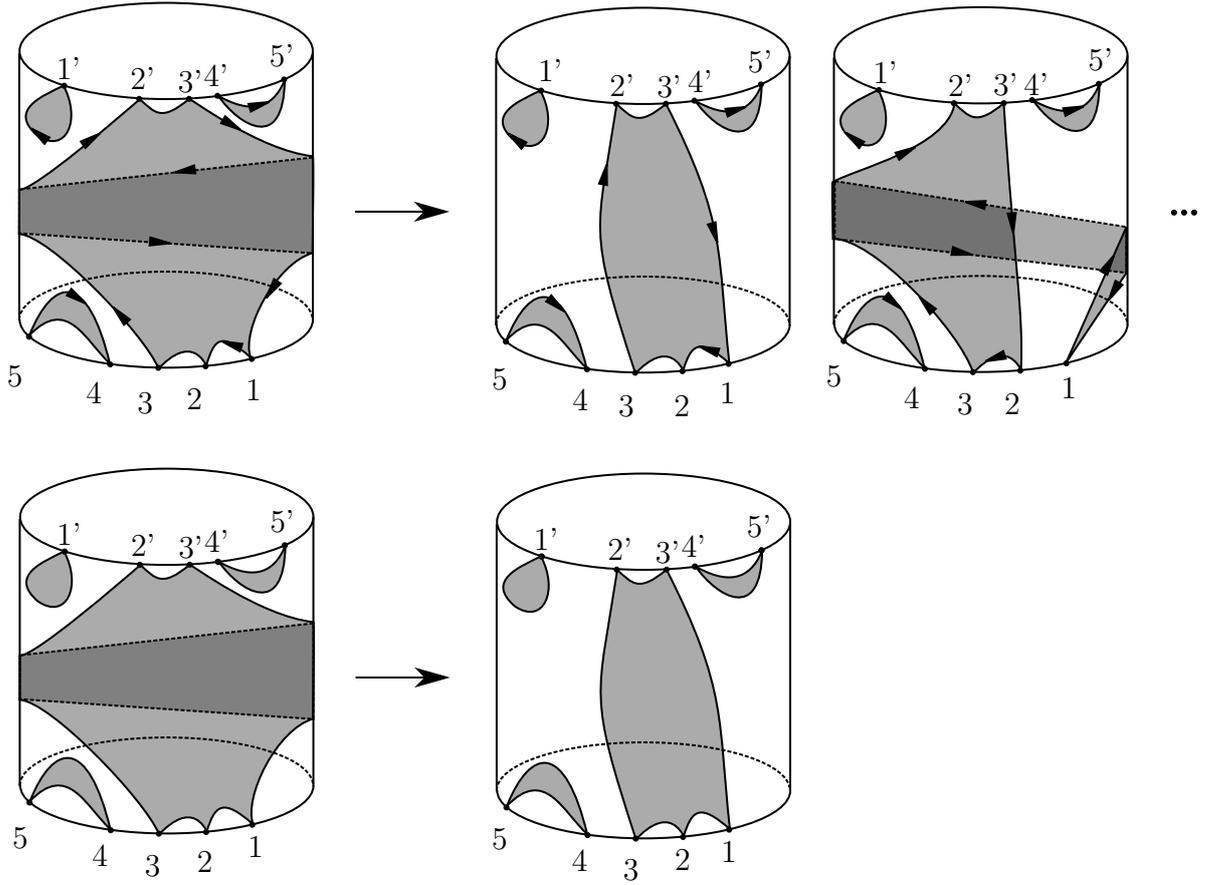}
\caption{The permutations above show how one can construct the case (3) from the case (1) by sending $\kappa_{i,j}\mapsto i\cdot j\cdot  \kappa_{i+j}$ since there are $i\cdot j$ possibilities to generate cycles of length $i+j$ where all the other trivial cycles are the same. The partitions below show that from the case (3) just one partition can be generated $\kappa_{i,j}\mapsto \kappa_{i+j}$. All $i\cdot j$ permutations above are in the same equivalence class as a partition.}
\label{Fig:cylinderPart}
\end{figure}

Therefore using the formula \eqref{W02Perm} we have to subtract the overcounted permutations to get partitions by setting
\begin{align*}
    \kappa_{i,j}\mapsto \kappa_{i,j}+(1-i\cdot j)\kappa_{i+j}.
\end{align*}
A straightforward computation yields
\begin{align*}
    &\sum_{i,j\geq 1}(1-i\cdot j)\kappa_{i+j} y_1^{i-1}y_2^{j-1}\\
    =&\sum_{n =2}^\infty \kappa_n \sum_{k=1}^{n-1}(1-k (n-k))y_1^{k-1}y_2^{n-k-1}\\
    =&\frac{1}{y_1y_2}(1-y_1\partial_{y_1}y_2\partial_{y_2}) \bigg(y_1y_2\frac{X(y_1)-X(y_2)}{y_1-y_2}+\frac{1}{y_1y_2}\bigg)\\
\end{align*}
which is the additional term appearing in \eqref{W02part} in comparison to \eqref{W02Perm}.

\end{proof}
\end{theorem}

The proof of Theorem \ref{Thm:W02Part} is somehow equivalent to the statements in \cite{10.1155/S1073792804133023}, but combined additionally with the functional relation \eqref{W02Perm} for permutations on the cylinder. However, the language used in \textit{op. cit.}
is much more technical but follows the same idea. The coefficients can be extracted easily from the generating series via
\begin{align}\label{w02m}
    m^{(0)}_{i,j}=\Res_{y_1,y_2\to 0} W^{(0)}_{par,2}(X(y_1),X(y_2)) X(y_1)^iX(y_2)^j X'(y_1)X'(y_2) dy_1dy_2.
\end{align}
\begin{example}
    The first coefficients are computed from \eqref{w02m}
    \begin{align*}
        m^{(0)}_{1,1}=&\kappa_{1,1}+\kappa_{2}\\
        m^{(0)}_{1,2}=&\kappa_{2,1}+2\kappa_{1}\kappa_{1,1}+\kappa_{3}+2\kappa_{1}\kappa_{2}\\
        m^{(0)}_{2,2}=&\kappa_{2,2}+4\kappa_{1}\kappa_{2,1}+4\kappa_{1}^2\kappa_{1,1}+\kappa_{4}+4\kappa_1\kappa_3+2\kappa_2^2+4\kappa_1^2\kappa_2\\
        m^{(0)}_{1,3}=&\kappa_{1,3}+3 \kappa_{1} \kappa_{1,2}+3 \kappa_{2} \kappa_{1,1}+3 \kappa_{1}^{2} \kappa_{1,1}+\kappa_{4}+3 \kappa_{1} \kappa_{3}+3 \kappa_{1}^{2} \kappa_{2}+3 \kappa_{2}^{2}\\
        m^{(0)}_{2,3}=&\kappa_{2,3}+3 \kappa_{1} \kappa_{2,2}+2 \kappa_{1} \kappa_{1,3}+3 \kappa_{2} \kappa_{1,2}+9 \kappa_{1}^{2} \kappa_{1,2}+6 \kappa_{1} \kappa_{2} \kappa_{1,1}+6 \kappa_{1}^{3} \kappa_{1,1}\\
        &+\kappa_{5}+5 \kappa_{1} \kappa_{4}+9 \kappa_{2} \kappa_{3}+9 \kappa_{1}^{2} \kappa_{3}+6 \kappa_{1}^{3} \kappa_{2}+12 \kappa_{1} \kappa_{2}^{2}\\
        m^{(0)}_{3,3}=&\kappa_{3,3}+6 \kappa_{1} \kappa_{3,2}+6 \kappa_{1}^{2} \kappa_{1,3}+6 \kappa_{2} \kappa_{1,3}+9 \kappa_{1}^{2} \kappa_{2,2}+18 \kappa_{1} \kappa_{2} \kappa_{1,2}+18 \kappa_{1}^{3} \kappa_{1,2}+9 \kappa_{2}^{2} \kappa_{1,1}\\
        &+18 \kappa_{1}^{2} \kappa_{2} \kappa_{1,1}+9 \kappa_{1}^{4} \kappa_{1,1}+\kappa_{6}+6 \kappa_{1} \kappa_{5}+15 \kappa_{2} \kappa_{4}
        +9 \kappa_{1}^{4} \kappa_{2}+18 \kappa_{1}^{3} \kappa_{3}+36 \kappa_{1}^{2} \kappa_{2}^{2}\\
        &+9 \kappa_{3}^{2}+15 \kappa_{1}^{2} \kappa_{4}+54 \kappa_{1} \kappa_{2} \kappa_{3}+12 \kappa_{2}^{3}
    \end{align*}
    which we also get from \cite[Remark 6.6]{Collins2006SecondOF} by substituting $\kappa_{i,j}\mapsto \kappa_{i,j}+(1-i\cdot j)\kappa_{i+j}$.
\end{example}
\begin{remark}
    Genus permutations and it generalisation with topologically nontrivial cycles $\kappa_{n_1,...,n_b}^{(g)}$ is known to be governed by Topological Recursion, see Remark \ref{rem:gbperm}. Topological Recursion as an algorithm works recursively in the Euler characteristic $2g+b-2$, where $b$ is the number of boundary components. It is very natural to ask if also genus partitions including more boundaries and more general topologically nontrivial blocks follow an analogous recursion in the Euler characteristic. The method which derived the genus 2 result in \cite{zuber2023counting} is not of this nature. It is worth to consider all generating series of partitions for any $g$ and $b$ and find relations between them, which are in this context called loop equations (Tutte equations or Dyson-Schwinger equations).
\end{remark}

\appendix
\section{Examples for Genus Permutations}\label{appper}
In the following, we want to list a few important examples by specialising $\kappa_i$. The inverse of $X(y)$ will have just two branches such that $W^{(g)}_{per}(x)$ can be written explicitly. 

\subsection{Factorials $\kappa_i=1$} 
These results are of course not new at all (see for instance \cite{MR314686,MR314687,MR0404045}), but achieved with new techniques and in a universal setting.

For $\kappa_i=1$, the geometric series yields
\begin{align*}
    X(y)=\frac{1}{y}+\frac{1}{1-y}.
\end{align*}
Inserting this into Theorem \ref{ThmPerm}, we find 
\begin{align}
    \mathcal{W}_{per}(X(y))=\sum_{m\geq 0} &\bigg(\frac{y^2(1-y)^2}{1-2y}\frac{\partial}{\partial {y}}\bigg)^m\bigg(\frac{y^2(1-y)^2}{1-2y}\bigg)\\\nonumber
    &\times
	[u^m]\frac{\exp\bigg(\frac{1}{\hbar}\log \bigg(\frac{(y+\frac{\hbar u}{2}) (y-1-\frac{\hbar u}{2})}{(y-\frac{\hbar u}{2})(y-1+\frac{\hbar u}{2})}\bigg)-\bigg(\frac{1}{y}+\frac{1}{1-y}\bigg)u\bigg)}{\hbar u},
\end{align}
where $\mathcal{W}_{per}=\sum_g\hbar^{2g-1}W^{(g)}_{1}$. This closed-type formula is to our knowledge new.

Substituting $X(y)=\frac{1}{y}+\frac{1}{1-y}$ into Example \ref{experg}, we get
\begin{align*}
    W^{(0)}_{per}(x)=&\frac{1}{2}-\frac{1}{2} \sqrt{\frac{x-4}{x}}\\
    =&\frac{1}{x}+\frac{1}{x^2}+\frac{2}{x^3}+\frac{5}{x^4}+\frac{14}{x^5}+\frac{42}{x^6}+\frac{132}{x^7}+\frac{429}{x^8}+\frac{1430}{x^9}+\frac{4862}{x^{10}}+\mathcal{O}(x^{-11})\\
    W^{(1)}_{per}(x)=&\frac{1}{(x-4)^{5/2} x^{3/2}}\\
    =&\frac{1}{x^4}+\frac{10}{x^5}+\frac{70}{x^6}+\frac{420}{x^7}+\frac{2310}{x^8}+\frac{12012}{x^9}+\frac{60060}{x^{10}}+\mathcal{O}(x^{-11})\\
    W^{(2)}_{per}(x)=&\frac{8 x^2-8 x+9}{(x-4)^{11/2} x^{5/2}}\\
    =&\frac{8}{x^6}+\frac{168}{x^7}+\frac{2121}{x^8}+\frac{20790}{x^9}+\frac{174174}{x^{10}}+\mathcal{O}(x^{-11})\\
    W^{(3)}_{per}(x)=&\frac{180 x^{4}-32 x^{3}+528 x^{2}-720 x+450}{(x-4)^{17/2} x^{7/2}}\\
    =&\frac{180}{x^{8}}+\frac{6088}{x^{9}}+\frac{115720}{x^{10}}+\mathcal{O}(x^{-11}).
\end{align*}
Adding them together, we find 
\begin{align*}
    \hbar \mathcal{W}_{per}(x)=&\frac{1}{x}+\frac{1}{x^2}+\frac{2}{x^3}+\frac{5+\hbar^2}{x^4}+\frac{14+\hbar^210}{x^5}+\frac{42+\hbar^270+\hbar^4 8}{x^6}\\
    +&\frac{132+\hbar^2420+\hbar^4 168}{x^7}+\frac{429+\hbar^2 2310+\hbar^4 2121+ \hbar^6 180}{x^8}\\
    +&\frac{1430+\hbar^2 12012+\hbar^4 20790+ \hbar^6 6088}{x^9}+\mathcal{O}(x^{-10})
\end{align*}
Note that at $\hbar=1$, this sums up to the factorials as a formal expression coefficient-wise
\begin{align*}
     \mathcal{W}_{per}(x)\vert_{\hbar=1}=W_{per}(x)=\sum_{n=0}^\infty \frac{n!}{x^{n+1}}.
\end{align*}

\subsection{Stirling Numbers of the first kind $\kappa_i=\kappa$}
Enhancing the problem to Stirling numbers of the first kind, we take $\kappa_i=\kappa$. The geometric series yields
\begin{align*}
    X(y)=\frac{1}{y}+\frac{\kappa}{1-y}.
\end{align*}
Inserting this into Theorem \ref{ThmPerm}, we find 
\begin{align}
    \mathcal{W}_{per}(X(y))=\sum_{m\geq 0} &\bigg(\frac{1}{\frac{1}{y^{2}}-\frac{\kappa}{\left(1-y\right)^{2}}}\frac{\partial}{\partial {y}}\bigg)^m\bigg(\frac{1}{\frac{1}{y^{2}}-\frac{\kappa}{\left(1-y\right)^{2}}}\bigg)\\\nonumber
    &\times
	[u^m]\frac{\exp\bigg(\frac{1}{\hbar}\log \bigg(\frac{(y+\frac{\hbar u}{2}) (y-1-\frac{\hbar u}{2})^\kappa}{(y-\frac{\hbar u}{2})(y-1+\frac{\hbar u}{2})^\kappa}\bigg)-\bigg(\frac{1}{y}+\frac{\kappa}{1-y}\bigg)u\bigg)}{\hbar u},
\end{align}
where $\mathcal{W}_{per}=\sum_g\hbar^{2g-1}W^{(g)}_{1}$. This closed-type formula is to our knowledge also new.

Inserting $X(y)=\frac{1}{y}+\frac{\kappa}{1-y}$ into Example \ref{experg}, we get
\begin{align*}
    W^{(0)}_{per}(x)=&\frac{-\kappa+x+1-\sqrt{\kappa^{2}-2\left(x+1\right) \kappa+\left(-1+x\right)^{2}}}{2 x}\\
    =&\frac{1}{x}+\frac{\kappa}{x^{2}}+\frac{\kappa \left(\kappa+1\right)}{x^{3}}+\frac{\kappa^{3}+3 \kappa^{2}+\kappa}{x^4}+\frac{\kappa \left(\kappa+1\right) \left(\kappa^{2}+5 \kappa+1\right)}{x^5}+\mathcal{O}(x^{-6})\\
    W^{(1)}_{per}(x)=&\frac{\kappa x}{\left(\kappa^{2}-2\left( x+1\right) \kappa+\left(-1+x\right)^{2}\right)^{5/2}}\\
    =&\frac{\kappa}{x^{4}}+\frac{5 \kappa \left(\kappa+1\right)}{x^5}+\frac{15 \kappa^{3}+40 \kappa^{2}+15 \kappa }{x^6}+\frac{35 \kappa  \left(\kappa +1\right) \left(\kappa^{2}+4 \kappa +1\right)}{x^7}+\mathcal{O}(x^{-8})\\
    W^{(2)}_{per}(x)=&\frac{\left(8 x^{4}-4\left(\kappa+1\right) x^{3}+\left(-15 (\kappa-1)^{2}+9 \kappa\right) x^{2}+10(\kappa+1)(\kappa-1)^2 x+(\kappa-1)^{4}\right) \kappa x}{\left(\kappa^{2}-2\left( x+1\right) \kappa+\left(-1+x\right)^{2}\right)^{11/2}}\\
    =&\frac{8 \kappa }{x^6}+\frac{84 \kappa^{2}+84 \kappa }{x^7}+\frac{469 \kappa^{3}+1183 \kappa^{2}+469 \kappa}{x^{8}}+\mathcal{O}(x^{-9})
\end{align*}
Adding them together, we find 
\begin{align*}
    \hbar \mathcal{W}_{per}(x)=&\frac{1}{x}+\frac{\kappa}{x^2}+\frac{\kappa \left(\kappa+1\right)}{x^3}+\frac{\kappa^{3}+3 \kappa^{2}+\kappa(1+\hbar^2) }{x^4}\\
    &+\frac{\kappa \left(\kappa+1\right) \left(\kappa^{2}+5 \kappa+1+5 \hbar^2\right)}{x^5}+\mathcal{O}(x^{-6})
\end{align*}
Note again that at $\hbar=1$, this sums up to the Stirling numbers of the first kind as a formal expression coefficient-wise
\begin{align*}
     \mathcal{W}_{per}(x)\vert_{\hbar=1}=W_{per}(x)=\sum_{n=0}^\infty \frac{\sum_{k=0}^n(-1)^{n-k}s(n,k) \kappa^k}{x^{n+1}}.
\end{align*}

\subsection{Harer-Zagier $\kappa_i=\delta_{i,2}$}
Specialising the problem to count the Euler characteristic of the moduli space of complex curves, we take $\kappa_i=\delta_{i,2}$, i.e.
\begin{align*}
    X(y)=\frac{1}{y}+y.
\end{align*}
Inserting this into Theorem \ref{ThmPerm}, we find 
\begin{align}\label{BPerm}
    \mathcal{W}_{per}(X(y))=\sum_{m\geq 0} &\bigg(\frac{y^2}{1-y^2}\frac{\partial}{\partial {y}}\bigg)^m\bigg(\frac{y^2}{1-y^2}\bigg)\\\nonumber
    &\times
	[u^m]\frac{\exp\bigg(\frac{1}{\hbar}\log \bigg(\frac{(y+\frac{\hbar u}{2})}{(y-\frac{\hbar u}{2})}\bigg)-\frac{u}{y}\bigg)}{\hbar u},
\end{align}
where $\mathcal{W}_{per}=\sum_g\hbar^{2g-1}W^{(g)}_{1}$. This closed-type formula is to our knowledge also new.

Inserting $X(y)=\frac{1}{y}+y$ into Example \ref{experg}, we get
\begin{align*}
    W^{(0)}_{per}(x)=&\frac{1}{2} \left(x-\sqrt{x^2-4}\right)\\
    =&\frac{1}{x}+\frac{1}{x^3}+\frac{2}{x^5}+\frac{5}{x^7}+\frac{14}{x^9}+\frac{42}{x^{11}}+\frac{132}{x^{13}}+\frac{429}{x^{15}}+\frac{1430}{x^{17}}+\frac{4862}{x^{19}}+\mathcal{O}(x^{-21})\\
    W^{(1)}_{per}(x)=&\frac{1}{\left(x^2-4\right)^{5/2}}\\
    =&\frac{1}{x^5}+\frac{10}{x^7}+\frac{70}{x^9}+\frac{420}{x^{11}}+\frac{2310}{x^{13}}+\frac{12012}{x^{15}}+\frac{60060}{x^{17}}+\frac{291720}{x^{19}}+\mathcal{O}(x^{-21})\\
    W^{(2)}_{per}(x)=&\frac{21 \left(x^2+1\right)}{\left(x^2-4\right)^{11/2}}\\
    =&\frac{21}{x^9}+\frac{483}{x^{11}}+\frac{6468}{x^{13}}+\frac{66066}{x^{15}}+\frac{570570}{x^{17}}+\frac{4390386}{x^{19}}+\mathcal{O}(x^{-21})\\
    W^{(3)}_{per}(x)=&\frac{11 \left(135 x^4+558 x^2+158\right)}{\left(x^2-4\right)^{17/2}}\\
    =&\frac{1485}{x^{13}}+\frac{56628}{x^{15}}+\frac{1169740}{x^{17}}+\frac{17454580}{x^{19}}+\mathcal{O}(x^{-21})
\end{align*}
Adding them together, we find 
\begin{align*}
    \hbar \mathcal{W}_{per}(x)=&\frac{1}{x}+\frac{1}{x^3}
    +\frac{2+\hbar^2}{x^5}+\frac{5+\hbar^2 10}{x^7}+\frac{14+\hbar^2 70+\hbar^4 21}{x^9}+\frac{42+\hbar^2 420+\hbar^4 483}{x^{11}}\\
    +&\frac{132+\hbar^2 2310+\hbar^4 6468+\hbar^6 1485}{x^{13}}+\frac{429+\hbar^2 12012+\hbar^4 66066+\hbar^6 56628}{x^{15}}\\
    +&\mathcal{O}(x^{-17}).
\end{align*}
Note again that at $\hbar=1$, this sums up to the double factorial as a formal expression coefficient-wise
\begin{align*}
     \mathcal{W}_{per}(x)\vert_{\hbar=1}=W_{per}(x)=\sum_{n=0}^\infty \frac{(2n-1)!!}{x^{2n+1}}.
\end{align*}

\section{Examples for Genus Partitions}\label{app:part}
Examples for genus partitions can be found extensively in \cite{coquereaux2023counting}. To connect their examples with our notation, we have to specify $X(y)$ in Proposition \ref{PropPart} by the following way:
\begin{align*}
    &\text{Bell numbers:} &&X(y)=\frac{1}{y}+\frac{1}{1-y}\\
    &\text{Stirling 2nd:} &&X(y)=\frac{1}{y}+\frac{\kappa}{1-y}\\
    &\text{Harer-Zagier:} &&X(y)=\frac{1}{y}+y
\end{align*}
or with "no singletons" just subtract 1 from $X(y)$ from the first two examples. Note that for the Harer-Zagier example, the formulae for genus permutations and partitions coincide since the we are just allowing for permutations to have cycles of length two which can give just one partition. However, the general formulae for genus permutations and partitions are very different.

Exactly in the same way as in Appendix \ref{appper}, if we set $\hbar=1$ the genus Bell numbers and Stirling numbers of the second kind sum up coefficient-wise to the ordinary Bell and Stirling numbers.\\

Now, we prove some minor open conjectures of \cite{coquereaux2023counting}.
The following Proposition proves and generalises an observation in \cite[Sec. 5.1]{coquereaux2023counting} about so-called Faa di Bruno coefficients 
\begin{proposition}
    Fix the size of blocks for partitions to be $p\in \mathbb{N}$. Then, we have for $k\in \mathbb{N}$
    \begin{align*}
        m^{(1)}_{pk}=\frac{(p-1)^2 p}{2}\sum_{l=0}^{k-2} \binom{pk}{l}\frac{(k-1-l)(k-l)(k+1-l)}{6} (p-1)^{k-2-l}.
    \end{align*}
    \begin{proof}
        Choose $X(y)=\frac{1}{y}+y^{p-1}$ since any block of the partition has just length $p$. This gives from Proposition \ref{PropPart} and \eqref{res2g}
\begin{align*}
     m^{(1)}_{pk}=&\mathrm{Res}_{y=0} \frac{(p-1)^2 p y^{2 p-1}}{2 \left((p-1) y^p-1\right)^4} \bigg(\frac{1}{y}+y^{p-1}\bigg)^{p k}\\
    =&\frac{(p-1)^2 p}{2}\mathrm{Res}_{y=0} \frac{1}{y^{(k-2)p+1}}\frac{ \bigg(1+y^{p}\bigg)^{p k}}{ \left(1-(p-1) y^p\right)^4} .
\end{align*}
Let us expand the right term
\begin{align*}
    \frac{ \bigg(1+y^{p}\bigg)^{p k}}{ \left(1-(p-1) y^p\right)^4}=\bigg(\sum_{i=0}^{pk}\binom{pk}{i}y^{i p}\bigg)\bigg(\sum_{j=0}^\infty \frac{(j+1)(j+2)(j+3)}{6} (p-1)^jy^{jp}\bigg).
\end{align*}
Now, only the $y^{(k-2)p}$-coefficient contributes to the residue, which proves the assertion.
\end{proof}
\end{proposition}
Similar computations can be performed for $g=2$, but it is much more involved and does not give more insight.

Another conjecture in \cite[eq. (100)]{coquereaux2023counting} is about genus 1 partitions with exactly 3 blocks of length $r,p,q$.
\begin{proposition}
    For pairwise different $r,p,q\in \mathbb{N}$, we have
    \begin{align*}
        m^{(1)}_{r+p+q}=\frac{n}{2}\bigg[&p(p-1)(q+r-2)+q(q-1)(p+r-2)+r(r-1)(q+p-2) \\
    &+8(p-1)(r-1)(q-1)\bigg].
    \end{align*}
    \begin{proof}
        It is sufficient to take $X(y)=\frac{1}{y}+y^{r-1}+y^{p-1}+y^{q-1}$. Note for $n=r+p+q$, we can write 
        \begin{align*}
            X(y)^{n}=\frac{1}{y^n}(1+n(y^{p}+y^q+y^r)+n(n-1)(y^{r+q}+y^{p+q}+y^{r+p})+n(n-1)(n-2)y^n+...
        \end{align*}
        as an expansion about $y=0$, where all the other terms will not contribute. Inserting into Proposition \ref{PropPart} and \eqref{res2g}, we compute first 
\begin{align*}
    &W^{(1)}_{par}(X(y))X'(y)\\
    =&\frac{\partial}{\partial y}\bigg(\frac{1}{4(-1+(p-1)y^{p}+(q-1)y^{q}+(r-1)y^{r})^2}
    +\frac{1}{6(-1+(p-1)y^{p}+(q-1)y^{q}+(r-1)y^{r})^3}\bigg)\\
    =&\frac{(p-1)py^{p-1}+(q-1)qy^{q-1}+(r-1)ry^{r-1}}{2(1-(p-1)y^{p}-(q-1)y^{q}-(r-1)y^{r})^3}-\frac{(p-1)py^{p-1}+(q-1)qy^{q-1}+(r-1)ry^{r-1}}{2(1-(p-1)y^{p}-(q-1)y^{q}-(r-1)y^{r})^4}.
\end{align*}
Now collecting everything and taking the residue yields:
\begin{align*}
    &m^{(1)}_{r+p+q}=-\mathrm{Res}_{y=0} W^{(1)}(X(y)) X(y)^{n} X'(y)\\
    &=-\frac{1}{2}\bigg[-n(p(p-1)(q+r-2)+q(q-1)(p+r-2)+r(r-1)(q+p-2)) \\
    &\qquad +(12-20)(p+q+r)(p-1)(r-1)(q-1)\bigg]
\end{align*}
which is the assertion and equivalent to the statement in \cite[eq. (100)]{coquereaux2023counting}.
    \end{proof}
\end{proposition}

\newcommand{\etalchar}[1]{$^{#1}$}


\begin{thebibliography}{BCGF{\etalchar{+}}21b}
\expandafter\ifx\csname url\endcsname\relax
  \def\url#1{\texttt{#1}}\fi
\expandafter\ifx\csname doi\endcsname\relax
  \def\doi#1{\burlalt{doi:#1}{http://dx.doi.org/#1}}\fi
\expandafter\ifx\csname urlprefix\endcsname\relax\def\urlprefix{URL }\fi
\expandafter\ifx\csname href\endcsname\relax
  \def\href#1#2{#2}\fi
\expandafter\ifx\csname burlalt\endcsname\relax
  \def\burlalt#1#2{\href{#2}{#1}}\fi

\bibitem[ABDB{\etalchar{+}}22]{Alexandrov:2022ydc}
A.~Alexandrov, B.~Bychkov, P.~Dunin-Barkowski, M.~Kazarian, and S.~Shadrin.
\newblock {A universal formula for the $x-y$ swap in topological recursion}.
\newblock 12 2022, \burlalt{2212.00320}{http://arxiv.org/abs/2212.00320}.

\bibitem[BCGF21a]{Borot:2021eif}
G.~Borot, S.~Charbonnier, and E.~Garcia-Failde.
\newblock {Topological recursion for fully simple maps from ciliated maps}.
\newblock 6 2021, \burlalt{2106.09002}{http://arxiv.org/abs/2106.09002}.

\bibitem[BCGF{\etalchar{+}}21b]{Borot:2021thu}
G.~Borot, S.~Charbonnier, E.~Garcia-Failde, F.~Leid, and Sergey.
\newblock {Analytic theory of higher order free cumulants}.
\newblock 12 2021, \burlalt{2112.12184}{http://arxiv.org/abs/2112.12184}.

\bibitem[BDBKS20a]{Bychkov:2020ujd}
B.~Bychkov, P.~Dunin-Barkowski, M.~Kazarian, and S.~Shadrin.
\newblock {Explicit closed algebraic formulas for Orlov-Scherbin $n$-point
  functions}.
\newblock 8 2020, \burlalt{2008.13123}{http://arxiv.org/abs/2008.13123}.

\bibitem[BDBKS20b]{Bychkov:2020yzy}
B.~Bychkov, P.~Dunin-Barkowski, M.~Kazarian, and S.~Shadrin.
\newblock {Topological recursion for Kadomtsev-Petviashvili tau functions of
  hypergeometric type}.
\newblock 12 2020, \burlalt{2012.14723}{http://arxiv.org/abs/2012.14723}.

\bibitem[BDBKS21]{Bychkov:2021hfh}
B.~Bychkov, P.~Dunin-Barkowski, M.~Kazarian, and S.~Shadrin.
\newblock {Generalised ordinary vs fully simple duality for $n$-point functions
  and a proof of the Borot--Garcia-Failde conjecture}.
\newblock 6 2021, \burlalt{2106.08368}{http://arxiv.org/abs/2106.08368}.

\bibitem[BIPZ78]{Brezin:1977sv}
E.~Brezin, C.~Itzykson, G.~Parisi, and J.~Zuber.
\newblock {Planar Diagrams}.
\newblock {\em Commun. Math. Phys.}, 59:35, 1978.
\newblock \doi{10.1007/BF01614153}.

\bibitem[BR24]{brestensky2024noncrossing}
L.~G. Brestensky and N.~Reading.
\newblock Noncrossing partitions of an annulus, 2024,
  \burlalt{2212.14151}{http://arxiv.org/abs/2212.14151}.

\bibitem[CFF13]{CHAPUY20132064}
G.~Chapuy, V.~Féray, and E.~Fusy.
\newblock A simple model of trees for unicellular maps.
\newblock {\em Journal of Combinatorial Theory, Series A}, 120(8):2064--2092,
  2013.
\newblock \doi{https://doi.org/10.1016/j.jcta.2013.08.003}.

\bibitem[CH13]{MR3261809}
R.~Cori and G.~Hetyei.
\newblock Counting genus one partitions and permutations.
\newblock {\em S\'{e}m. Lothar. Combin.}, 70:Art. B70e, 29, 2013,
  \burlalt{1306.4628}{http://arxiv.org/abs/1306.4628}.

\bibitem[CH18]{MR3891093}
R.~Cori and G.~Hetyei.
\newblock Counting partitions of a fixed genus.
\newblock {\em Electron. J. Combin.}, 25(4):Paper No. 4.26, 37, 2018.
\newblock \doi{10.37236/7632}.

\bibitem[Cha11]{CHAPUY2011874}
G.~Chapuy.
\newblock A new combinatorial identity for unicellular maps, via a direct
  bijective approach.
\newblock {\em Advances in Applied Mathematics}, 47(4):874--893, 2011.
\newblock \doi{https://doi.org/10.1016/j.aam.2011.04.004}.

\bibitem[CJ10]{MR2550161}
S.~Cautis and D.~M. Jackson.
\newblock On {T}utte's chromatic invariant.
\newblock {\em Trans. Amer. Math. Soc.}, 362(1):491--507, 2010.
\newblock \doi{10.1090/S0002-9947-09-04836-3}.

\bibitem[CMSS07]{Collins2006SecondOF}
B.~Collins, J.~A. Mingo, P.~Sniady, and R.~Speicher.
\newblock Second order freeness and fluctuations of random matrices, iii.
  higher order freeness and free cumulants.
\newblock {\em Doc. Math.}, 12:1--70, 2007.

\bibitem[Cor75]{MR0404045}
R.~Cori.
\newblock {\em Un code pour les graphes planaires et ses applications}.
\newblock Ast\'{e}risque, No. 27. Soci\'{e}t\'{e} Math\'{e}matique de France,
  Paris, 1975.
\newblock With an English abstract.

\bibitem[Cvi81]{Cvitanovic:1980jz}
P.~Cvitanovic.
\newblock {Planar perturbation expansion}.
\newblock {\em Phys. Lett. B}, 99:49--52, 1981.
\newblock \doi{10.1016/0370-2693(81)90801-7}.

\bibitem[CZ23]{coquereaux2023counting}
R.~Coquereaux and J.-B. Zuber.
\newblock Counting partitions by genus. ii. a compendium of results, 2023,
  \burlalt{2305.01100}{http://arxiv.org/abs/2305.01100}.

\bibitem[EO07]{Eynard:2007kz}
B.~Eynard and N.~Orantin.
\newblock {Invariants of algebraic curves and topological expansion}.
\newblock {\em Commun. Num. Theor. Phys.}, 1:347--452, 2007,
  \burlalt{math-ph/0702045}{http://arxiv.org/abs/math-ph/0702045}.
\newblock \doi{10.4310/CNTP.2007.v1.n2.a4}.

\bibitem[GS98]{MR1649966}
A.~Goupil and G.~Schaeffer.
\newblock Factoring {$n$}-cycles and counting maps of given genus.
\newblock {\em European J. Combin.}, 19(7):819--834, 1998.
\newblock \doi{10.1006/eujc.1998.0215}.

\bibitem[Hoc22a]{Hock:2022pbw}
A.~Hock.
\newblock {A simple formula for the $x$-$y$ symplectic transformation in
  topological recursion}.
\newblock 11 2022, \burlalt{2211.08917}{http://arxiv.org/abs/2211.08917}.

\bibitem[Hoc22b]{Hock:2022wer}
A.~Hock.
\newblock {On the $x$-$y$ Symmetry of Correlators in Topological Recursion via
  Loop Insertion Operator}.
\newblock 1 2022, \burlalt{2201.05357}{http://arxiv.org/abs/2201.05357}.

\bibitem[Hoc23]{Hock:2023qii}
A.~Hock.
\newblock {Laplace transform of the $x-y$ symplectic transformation formula in
  Topological Recursion}.
\newblock 4 2023, \burlalt{2304.03032}{http://arxiv.org/abs/2304.03032}.

\bibitem[HZ86]{MR848681}
J.~Harer and D.~Zagier.
\newblock The {E}uler characteristic of the moduli space of curves.
\newblock {\em Invent. Math.}, 85(3):457--485, 1986.
\newblock \doi{10.1007/BF01390325}.

\bibitem[Kre72]{MR309747}
G.~Kreweras.
\newblock Sur les partitions non crois\'{e}es d'un cycle.
\newblock {\em Discrete Math.}, 1(4):333--350, 1972.
\newblock \doi{10.1016/0012-365X(72)90041-6}.

\bibitem[KZ15]{Kazarian:2014ooa}
M.~Kazarian and P.~Zograf.
\newblock {Virasoro constraints and topological recursion for Grothendieck's
  dessin counting}.
\newblock {\em Lett. Math. Phys.}, 105(8):1057--1084, 2015,
  \burlalt{1406.5976}{http://arxiv.org/abs/1406.5976}.
\newblock \doi{10.1007/s11005-015-0771-0}.

\bibitem[LZ04]{MR2036721}
S.~K. Lando and A.~K. Zvonkin.
\newblock {\em Graphs on surfaces and their applications}, volume 141 of {\em
  Encyclopaedia of Mathematical Sciences}.
\newblock Springer-Verlag, Berlin, 2004.
\newblock \doi{10.1007/978-3-540-38361-1}.
\newblock With an appendix by Don B. Zagier, Low-Dimensional Topology, II.

\bibitem[MN04]{10.1155/S1073792804133023}
J.~A. Mingo and A.~Nica.
\newblock {Annular noncrossing permutations and partitions, and second-order
  asymptotics for random matrices}.
\newblock {\em International Mathematics Research Notices},
  2004(28):1413--1460, 01 2004.
\newblock \doi{10.1155/S1073792804133023}.

\bibitem[Voi86]{Voiculescu1986AdditionOC}
D.~Voiculescu.
\newblock Addition of certain non-commuting random variables.
\newblock {\em Journal of Functional Analysis}, 66:323--346, 1986.

\bibitem[WL72a]{MR314687}
T.~Walsh and A.~B. Lehman.
\newblock Counting rooted maps by genus. {II}.
\newblock {\em J. Combinatorial Theory Ser. B}, 13:122--141, 1972.
\newblock \doi{10.1016/0095-8956(72)90049-4}.

\bibitem[WL72b]{MR314686}
T.~R.~S. Walsh and A.~B. Lehman.
\newblock Counting rooted maps by genus. {I}.
\newblock {\em J. Combinatorial Theory Ser. B}, 13:192--218, 1972.
\newblock \doi{10.1016/0095-8956(72)90056-1}.

\bibitem[Zub23]{zuber2023counting}
J.-B. Zuber.
\newblock Counting partitions by genus. i. genus 0 to 2, 2023,
  \burlalt{2303.05875}{http://arxiv.org/abs/2303.05875}.

\end{thebibliography}
\end{document}